\documentclass[10pt,a4paper]{amsart}
\usepackage{amssymb}
\usepackage{amsmath}
\usepackage{graphicx}
\usepackage{unbm}

\begin{document}
\title{Some properties of the symbol algebras}
\author{Diana Savin, Cristina Flaut, Camelia Ciobanu}
\address{
University  "Ovidius"\newline
\indent Department of Mathematics and Informatics \newline
\indent Bd. Mamaia 124, 900527, Constanta, Romania}
\email{savin.diana@univ-ovidius.ro, cflaut@univ-ovidius.ro, cristina\_flaut@yahoo.com}
\address{Department of Mathematics-Informatics and Fundamental
Technical Sciences\newline 
\indent Mircea cel Batran Naval Academy\newline 
\indent 1, Fulgerului\ \  Sreet, 900218, Constanta, 
Romania}
\email{c\_cami\_ro@yahoo.com}
\setcounter{page}{1}\coordinates{0}{0}{00 - 00}

\keywords{\em symbol algebras, p-adic valuation, Artin symbol\\
Accepted on June 2, 2009}
{\date{7.11. 2008}\undate{1}\adate{1}}
\footnote{For the first author, the research for this article was carried out during her visit at the Central European University (CEU), Budapest, from 1 January to 31 March 2008, and supported by the CEU Special and Extension Program.}
\maketitle

\begin{unabstract}In this paper, we obtain some properties of the
symbol algebras, starting  from their
connections with the quaternion and cyclic algebras over a field $K_{p},$%
where $K$ is an algebraic number field, $p$ is
a prime in $K$ and $K_{p}$ is the
completion of  $K$   with 
respect to $p-$  adic valuation,  in the case when
$K=\mathbb{Q}\left( \varepsilon \right) ,\,\,\varepsilon ^{3}=1,\varepsilon
\neq 1.$
\end{unabstract}

\[
\]

{\large 1. Introduction} 
\[
\]

Symbol algebras have many applications in number theory (class field theory), as can be seen in [4], [6], [7]. Since they are a  natural generalization of the quaternion algebras, in this paper we find some interesting example of split quaternion algebras and non division symbol algebras and we give a necessary and sufficient condition for \thinspace a $%
K_{v}-$cyclic central simple algebra $A=\left( \frac{\alpha ,\beta }{%
K,\epsilon }\right) $ to be a division algebra. 

First, we recall some definitions in the theory of associative algebras.

Let $A\neq 0$ be an algebra over the field $K.$ If the equations $%
ax=b,\,ya=b,\forall a,b\in A,a\neq 0,$ have unique solutions, then the
algebra $A$ is called a \textit{division algebra. }If $A$ is a
finite-dimensional algebra, then $A$ is a division algebra if and only if $A$
is without zero divisors\thinspace ($x\neq 0,y\neq 0\Rightarrow xy\neq 0$).(see [9])

Let $K$ be a field with $charK\neq 2.$ Let $\mathbb{H}_{K}\left( \alpha ,\beta
\right) $ be a quaternion algebra with basis $\{1,e_{1},e_{2},e_{3}\}\,$%
\thinspace and the multiplication given by

{\footnotesize 
\[
\begin{tabular}{l|llll}
$\cdot $ & $1$ & $\,\,\,e_{1}$ & $\,\,e_{2}$ & $\,\,\,\,e_{3}$ \\ \hline
$1$ & $1$ & $e_{1}$ & $e_{2}$ & $\,\,\,e_{3}$ \\ 
$e_{1}$ & $e_{1}$ & $\alpha $ & $e_{3}$ & $\alpha e_{2}$ \\ 
$e_{2}$ & $e_{2}$ & $-e_{3}$ & $\beta $ & $-\,\,\,\beta e_{1}$ \\ 
$e_{3}$ & $e_{3}$ & $-\alpha e_{2}$ & $\beta e_{1}$ & $-\alpha \beta $%
\end{tabular}
. 
\]
}Each element $x\in \mathbb{H}_{K}\left( \alpha ,\beta \right) \,\,$has the
form $x=x_{0}\cdot 1+x_{1}e_{2}+x_{2}e_{2}+x_{3}e_{3},\,$with $x_{i}\in
K,\,i=0,1,2,3.\,\,\,\,$For $a\in \mathbb{H}_{K}\left( \alpha ,\beta \right)
,\,\,a=\,a_{0}+a_{1}e_{1}+a_{2}e_{2}+a_{3}e_{3},\,$ the element $\bar{a}%
=a_{0}-a_{1}e_{1}-a_{2}e_{2}-a_{3}e_{3}$ is called the \textit{conjugate} of
the element $a.$ Let $a\in \mathbb{H}_{K}\left( \alpha ,\beta \right) .$ We
have that $t\left( a\right) \cdot 1=a+\overline{a}\in K,\,n\left( a\right)
\cdot 1=a\overline{a}\in K$ and these are called the \textit{trace},
respectively, the \textit{norm} of the element $a\in $ $A$ .\thinspace
\thinspace \thinspace It\thinspace \thinspace \thinspace follows\thinspace
\thinspace \thinspace that$\,\,\left( a+\overline{a}\right) a\,=a^{2}+%
\overline{a}a$=$a^{2}+n\left( a\right) \cdot 1$ and\thinspace \thinspace $%
a^{2}-t\left( a\right) a+n\left( a\right) =0,\forall a\in A,$ $\,$therefore
the generalized quaternion algebras are \textit{quadratic}. We remark that $%
n\left( a\right) =a_{0}^{2}-\alpha a_{1}^{2}-\beta a_{2}^{2}+\alpha \beta
a_{3}^{2}$ . The generalized quaternion algebras is a division algebra if
and only if for $x\in \mathbb{H}_{K}\left( \alpha ,\beta \right) $ we have $%
n\left( x\right) =0$ only for $x=0.\,$ Otherwise, the algebra $\mathbb{H}%
_{K}\left( \alpha ,\beta \right) $ is \textit{a split }algebra.

An important invariant for a quaternion algebra $\mathbb{H}_{K}\left( \alpha
,\beta \right) $ is the \textit{associated conic}, denoted $C\left( \alpha
,\beta \right) .$ The associated conic is the projective plane curve defined
by the homogeneous equation $\alpha x^{2}+\beta y^{2}=z^{2}.\medskip $

Let $K$ be an algebraic number field. By a prime of $K$ we mean a class of
equivalent valuations of $K.$ Recall that the finite primes of $K$ are in
one-to-one correspondence with the primes ideals of the ring of integers of $%
K,$ and the infinite primes are in correspondence with the embedding of $K$
into the field of complex numbers $\mathbb{C}.$ If $\,v$ is a prime of $K,$ we
denote with $K_{v}$ the completion of $K$ with respect to the v-adic
valuation.\medskip 

\textbf{Proposition 1.1.} [4, pag. 7]\textit{The quaternion algebra }$\mathbb{H}%
_{K}\left( \alpha ,\beta \right) $\textit{\ is split if and only if the
conic }$C\left( \alpha ,\beta \right) $\textit{\ has a rational points over }%
$K($\textit{i.e. if there are }$x_{0},y_{0},z_{0}\in K$\textit{\ such that }$%
\alpha x_{0}^{2}+\beta y_{0}^{2}=z_{0}^{2}).\medskip $

A natural generalization of \thinspace \thinspace the \textit{quaternion
algebra } is the \textit{symbol algebra}, also known as a \textit{power norm
residue algebra}. J. Milnor, in his book \textit{Introduction to Algebraic
K-Theory}, calls it \thinspace the \textit{symbol algebra} because of \thinspace its
connection with the $K-$theory and with the Steinberg symbols.(see [8])

A \textit{symbol algebra} is a unitary associative algebra over a field $K$
with $\zeta \in K,$ $\zeta ^{n}=1,\zeta $ a primitive root, generated by the
elements $x,y$ \thinspace which satisfy the relations $x^{n}=\alpha
,y^{n}=\beta $ and $yx=\zeta xy.$ This algebra is denoted $\left( \frac{\alpha ,\beta }{K,\zeta }\right) .\,\,$%

Obviously,\: for $n=2$ we obtain the algebra $\mathbb{H}_{K}\left( \alpha ,\beta
\right) .$

The quaternion generalized algebras and symbol algebras are central simple
algebras.$\medskip $

\textbf{Proposition 1.2.} [8, pag. 237] \textit{If }$K$\textit{\ is an
algebraic number field and }$A$\textit{\ is a central simple }$K-$\textit{%
algebra, then the dimension of }$A$\textit{\ over }$K$\textit{\ is a square}%
.\medskip 

\textbf{Definition 1.3.} Let $A$ be a central simple algebra of finite
dimension $n\,\,\,$over $K.$ The positive integer $d=\sqrt{n}$ is called 
\textit{the degree} of the algebra $A.\medskip $

\textbf{Theorem}(\textbf{Weddeburn}).[8, pag. 50] \textit{Let }$A$\textit{\
be a central simple algebra over the field}$\,\,K.$\textit{\ There are }$%
n\in \mathbb{N}^{*}$\textit{\ and a division algebra }$D,$\textit{\ }$%
K\subseteq D,$\textit{\ such that }$A\simeq $\textit{\ }$\mathcal{M}%
_{n}\left( D\right) .$\textit{\ The division algebra }$D$\textit{\ is unique
up to an isomorphism.\medskip }

\textbf{Definition 1.4.} With the notation of the above Theorem, the degree
of \thinspace the algebra $D$ over $K$ (as an algebra) is called \textit{the
index} of the algebra $A.\medskip $

For some $h\in N^{*},$ the tensor product over the field $K$ $A\otimes ...\otimes A$ ($h-$
times) is isomorphic to a full matrix algebra over $K.\medskip $

\textbf{Definition 1.5.} The smallest such an $h$ is called \textit{the
exponent} of the algebra $A.\medskip $

\textbf{Theorem 1.6.[}1\textbf{]} \textit{The algebra }$A$\textit{\ is a
division algebra if and only if its index and its degree are the
same.\medskip }

\textbf{Theorem 1.7.} \textbf{(Brauer-Hasse-Noether).[}8\textbf{]} \textit{%
Every central simple algebra over an algebraic number field is cyclic and
its index is equal to its exponent.\medskip } We shall use in the third
section some results from the theory of algebraic number fields and we
recall these here.\medskip \newline
\textbf{Theorem 1.8.} ([1]) \textit{Let }$K\subseteq E$\textit{\ be a cyclic
extension of commutative fields of degree }$d.$\textit{\ The cyclic }$K-$%
\textit{algebra }$A=\left( \frac{\alpha ,\beta }{K,\zeta }\right) $\textit{\
has the exponent }$d$\textit{\ if and only if }$\alpha \notin N_{L/K}\left(
L^{*}\right) ,$\textit{\ for each minimal subfield }$L$\textit{\ of }$E$%
\textit{\ over }$K.\medskip $

\textbf{Theorem 1.9.}([4]) \textit{Let }$K$\textit{\ be a field such that }$%
\zeta \in K,\,\,\zeta ^{n}=1,\zeta $\textit{\ is a primitive root, and let }$%
\alpha ,\beta \in K^{*}.$\textit{\ Then the following statements are
equivalent:}

\textit{i) The cyclic algebra }$A=\left( \frac{\alpha ,\beta }{K,\zeta }%
\right) $\textit{\ is split.}

\textit{ii) The element }$\beta $\textit{\ is a norm from the extension }$K
\subseteq K(\sqrt[n]{\alpha}) $ \medskip

\textbf{Theorem 1.10.} ([1;2;6]) \textit{Let }$K$\textit{\ be an algebraic
number field, }$v$\textit{\ be a prime of }$K$\textit{\ and }$K\subseteq L$%
\textit{\ a Galois extension. Let }$w$\textit{\ be a prime of }$L\,\,\,\,$%
\textit{lying above }$v$ \textit{\ such that }$K_{v}\subseteq L_{w}$\textit{%
\ is a unramified extension of }$K_{v}$\textit{\ of (residual) degree }$f.$%
\textit{\ Let }$b=\pi _{v}^{m}\cdot u_{v}\in K_{v}^{*},$ \textit{\ where } $%
\pi_{v}$ \textit{denote a prime element in}
 $ K_{v}$ \textit{and} $u_{v}$\textit{\ a unit in the ring of integers }$\mathcal{O}_{v},m\in \mathbb{Z}%
.$\textit{\ Then }$b\in N_{L_{w}\,/\,K_{v}}\left( L_{w}^{*}\right) $\textit{%
\ if and only if }$f\,\mid \,m.$\textit{\ In particular, every unit of }$%
\mathcal{O}_{v}$\textit{\ is the norm of a unit in }$L_{w}.\medskip $

\textbf{Theorem 1.11.} ([2;7]) \textit{Let }$K$\textit{\ be an algebraic
number field,}$\,\,e$\textit{\ be an admissible modulus of }$K$\textit{, }$v$%
\textit{\ be a finite prime of }$K,$\textit{\ }$v$\textit{\ divides }$e.$%
\textit{\ Let }$K\subseteq L$\textit{\ be a Galois extension. Let }$w$%
\textit{\ be any prime of }$L\,\,\,\,$\textit{lying above }$v.$\textit{\
Then an element }$a\in N_{L_{w}\,/\,K_{v}}\left( L_{w}^{*}\right) $\textit{\
if and only if the Artin symbol }$\left( \frac{%
\begin{array}{c}
\\ 
L_{w}\,/\,K_{v}
\end{array}
}{%
\begin{array}{c}
(a)
\end{array}
}\right) $\textit{\ is the identity in the Galois group }$Gal\left(
L_{w}\,/\,K_{v}\right),$ \textit{where} $(a)$ \textit{denotes the ideal generates by} $a$ \textit{in the ring of integers} $\mathcal{O}_{v} .\medskip$

\textbf{Theorem 1.12.} ([6]) \textit{Let }$\zeta $\textit{\ be a primitive
root of the unity of }$l-$\textit{order, where }$l$\textit{\ is a prime
natural number and let} $A$ \textit{be the ring of integers of the Kummer
field} $Q(\zeta, \sqrt[l]{\mu})$ . \textit{A prime ideal} $P$\textit{\ in
the ring} $\mathbb{Z}[\zeta ]$\textit{\ is in }$A$\textit{\ in one of the
situations:}

\textit{i) It is equal with the }$l-$\textit{power of a prime ideal from} $A,
$ \textit{if the} $l-$\textit{power character }$\left( \frac{\mu }{P}\right)
_{l}=0;$

\textit{ii) It is a prime ideal in} $A$, \textit{if} $\left( \frac{\mu }{P}%
\right) _{l}=$\textit{\ a rot of order }$l$\textit{\ of unity, different
from }$1$. \newline
\textit{iii) It decomposes in $l$ different prime ideals from} $A$, \textit{%
if} $\left( \frac{\mu }{P}\right) _{l}=1.$\medskip

\textbf{Theorem 1.13.} ([5;6]) \textit{Let }$l$\textit{\ be a natural number, }$%
l\geq 3$\textit{\ and }$\zeta $\textit{\ be a primitive root of the unity
of\thinspace \thinspace }$l$-\textit{order. If }$p$\textit{\ is a prime
natural number, }$l$\textit{\ is not divisible with }$p$\textit{\ and }$f$%
\textit{\ is the smallest positive integer such that }$p^{f}\equiv 1$\textit{%
\ mod }$l$\textit{, then we have} 
\[
p\mathbb{Z}[\zeta ]=P_{1}P_{2}....P_{r}, 
\]
\textit{where }$r=\frac{%
\begin{array}{c}
\\ 
\varphi \left( l\right)
\end{array}
}{%
\begin{array}{c}
f
\end{array}
},\varphi $\textit{\ is the Euler's function and }$P_{j},\,j=1,...,r$\textit{%
\ are different prime ideals in the ring }$\mathbb{Z}[\zeta ].$%
\[
\]

In the following, we consider the symbol algebra for $n=3$ and $K=\mathbb{Q}%
\left( \varepsilon \right) $ or $\mathbb{Q}_{p}\left( \varepsilon \right) ,$
where $\varepsilon $ is a primitive cubic root of unity and $p$ a prime
number.\medskip

\[
\]

{\large 2. Some example of quaternion and symbol algebras} 
\[
\]

\textbf{Proposition 2.1.} \textit{For }$\alpha =-1,\beta =p,p=4k+3,$\textit{%
\ a prime number, }$K=\mathbb{Q},$\textit{\ the algebra }$\mathbb{H}_{\mathbb{Q}%
}\left( -1,p\right) $\textit{\ is a division algebra.\medskip }

\textbf{Proof.} Let $\,\,x\in \mathbb{H}_{\mathbb{Q}}\left( -1,p\right)
,\,\,x=x_{0}\cdot 1+x_{1}e_{2}+x_{2}e_{2}+x_{3}e_{3},\,x_{i}\in \mathbb{Q}%
,\,i=0,1,2,3$ such that $n\left( x\right) =0.$ It results $%
x_{0}^{2}+x_{1}^{2}-px_{2}^{2}-px_{3}^{2}=0,$ then $p\mid
(x_{0}^{2}+x_{1}^{2}).$ Since $p=4k+3$ is a prime and\thinspace $p\mid
(x_{0}^{2}+x_{1}^{2}),$ we obtain that $p\mid (x_{2}^{2}+x_{3}^{2}),$ and
the powers of \thinspace \thinspace $p$ in the factorization of $%
\,\,x_{0}^{2}+x_{1}^{2}$ and $x_{2}^{2}+x_{3}^{2}$ are even. We obtain a
contradiction, therefore $x=0.\Box \medskip $

\textbf{Theorem (Gauss). }If $\,\,p\equiv 1$ \textit{mod} $3,$ then there
are integers $a,b$ such that $4p=a^{2}+27b^{2}.\medskip $

\textbf{Proposition 2.2.} \textit{If }$K=$\textit{\ }$\mathbb{Q}\left( \sqrt{3}%
\right) ,$\textit{\ then the quaternion algebra }$\mathbb{H}_{K}\left(
-1,p\right) ,$\textit{\ where }$p\equiv 1$\textit{\ mod }$3$\textit{\ is a
split algebra.\medskip }

\textbf{Proof.} Indeed, $\mathbb{H}_{K}\left( -1,p\right) $ is a split algebra
if and only if the associated conic $-x^{2}+py^{2}=z^{2}$ has $\mathbb{Q}\left( 
\sqrt{3}\right) -$ rational points. Using the Gauss's theorem, there are $%
a,b\in \mathbb{Z}$ such that $4p=a^{2}+27b^{2}.$ Then$\,\,$for$%
\,\,\,\,\,y_{0}=1,z_{0}=\frac{a}{2},\,x_{0}=\frac{3\sqrt{3}b}{2},$ the point 
$\left( \frac{3\sqrt{3}b}{2},1,\frac{a}{2}\right) $ is a $\mathbb{Q}\left( 
\sqrt{3}\right) -$ rational point for the associated \thinspace conic, and
we use Proposition 1.1.$\Box \medskip $

From the Wedderburn theorem, we know that a finite dimensional simple
algebra $A$ over a field $K$ is isomorphic to a matrix algebra $\mathcal{M}%
_{n}\left( D\right), $ \: for $D$ a division algebra. Let $K=\mathbb{Q}\left( \varepsilon \right)$ where $%
\varepsilon $ is a cubic root of unity  and let $d=[D:K]$ be the index of the algebra $A.$ The algebra 
$%
A=\left( \frac{\alpha ,\beta }{K,\zeta }\right) $ is a central simple
algebra of degree $3,$ hence $d \mid 3.$
 
For $\alpha =-1, \beta =1, $ the algebra $A$ is generated, for example, by the elements

 $X=\left( 
\begin{array}{lll}
-1 & 0 & 0 \\ 
0 & -\varepsilon  & 0 \\ 
0 & 0 & -\varepsilon ^{2}
\end{array}
\right) $ and $Y=\left( 
\begin{array}{lll}
0 & 1 & 0 \\ 
0 & 0 & 1 \\ 
1 & 0 & 0
\end{array}
\right) ,$\thinspace \thinspace \\
where $X^{3}=-1I_{3},$ $Y^{3}=I_{3}$ and $%
YX=\varepsilon XY.$ (see[3])  We obtain that  $A\simeq \mathcal{M}_{n}\left( \mathbb{Q}%
\left( \varepsilon \right) \right) .$  Therefore \thinspace \thinspace $d=1
$ and the algebra $A$ is not a division algebra.\medskip 

 We obtain the following
proposition\medskip 

\textbf{Proposition 2.3. } \textit{The algebras }$A=\left( \frac{\alpha
,\beta }{\mathbb{Q}\left( \varepsilon \right) ,\varepsilon }\right) ,$\textit{\
for }$\alpha ,\beta \in \{-1,1\}$ \textit{are not division algebras}%
.\medskip \textit{\ }

\textbf{Proof.} The algebra $A$ has dimension $9,$ hence degree $3,$ with
basis \\$B=\{1,x,y,x^{2},y^{2},xy^{2},xy,x^{2}y,x^{2}y^{2}\},$ $%
x^{3}=a,y^{3}=b.$ With the correspondence $x\rightarrow X,y\rightarrow Y,$
we have that $A\simeq \mathcal{M}_{n}\left( \mathbb{Q}\left( \varepsilon
\right) \right) ,$ the index $d=1\neq 3,$ where $3$ is the algebra's degree,
then $A$ is not a division algebra.(We used Proposition 1.6. )$\Box
\medskip $

If the central simple algebra $A$ is a division algebra, since has the
degree three, it results that it is a cyclic algebra. It results that there
are the elements $x\in A-K,$ $\alpha \in K\,\,\,$such that $x^{3}=\alpha \in
K.$ From the Noether-Skolem theorem, it results that there is an element $%
y\in A-K$ such that $yxy^{-1}=\varepsilon x.$ $\;$We have $y^{3}x=xy^{3}$
and $y^{3}y=yy^{3},$ then $y^{3}$ commutes with the generators $x,y,$
therefore $y^{3}\in K=C(A),$ the centralizer of the algebra $A.$ Hence,
there is $\beta \in K$ such that $y^{3}=\beta ,$ and $A=\left( \frac{\alpha
,\beta }{K,\zeta }\right) \simeq \mathcal{M}_{n}\left( D\right) ,$ with $%
[D:K]=3.$%
\[
\]

{\large 3. The algebra }$A=\left( \frac{\alpha ,\beta }{K_{v},\varepsilon }%
\right)$ \newline
\bigskip\newline
We consider the case of the algebra $A=\left( \frac{\alpha ,\beta }{%
K_{v},\varepsilon }\right)$ where $\varepsilon $ is a primitive cubic root
of unity. We give a necessary and sufficient condition for \thinspace a $%
K_{v}-$cyclic central simple algebra $A=\left( \frac{\alpha ,\beta }{%
K_{v},\varepsilon }\right) $ to be a division algebra and finally we find when $%
\beta $ is a norm for the field $K_{v}\left(\sqrt[3]{\alpha}\right),$
 where $K_{v}$ is the completion of the field $K $ with respect
the v-adic valuation.

Let $K$ be an algebraic number field and $v$ be a prime (finite or infinite)
of $\,K$ such that $\varepsilon \in K_{v},$ where $\varepsilon $ is a
primitive cubic root. We consider the $K_{v}-$ central simple algebra $%
A=\left( \frac{%
\begin{array}{c}
\\ 
\alpha ,\beta
\end{array}
}{%
\begin{array}{c}
K_{v},\varepsilon
\end{array}
}\right) ,$ $\alpha ,\beta \in K_{v}^{*}.\medskip $

\textbf{Proposition 3.1.} \textit{With the above notation, if\thinspace
\thinspace }$L=K\left(\sqrt[3]{\alpha}\right),$ \textit{\ the
following statement are equivalent:}

\textit{i) The algebra }$A=\left( \frac{%
\begin{array}{c}
\\ 
\alpha ,\beta
\end{array}
}{%
\begin{array}{c}
K_{v},\varepsilon
\end{array}
}\right) $\textit{\ is a division algebra.}

\textit{ii) }$\beta \notin N_{L_{w\,}/\,K_{v}}\left( L_{w}^{*}\right) ,$%
\textit{\ for each }$w\,$\textit{a prime of }$L$\textit{\ lying above }$%
v.$\newline
\medskip\newline
\textbf{Proof.} We consider the cyclic extension of fields $K_{v}\subseteq
L_{w}$ and we apply the Theorems 1.6, 1.7, 1.8. We obtain that the $K_{v}- $%
cyclic central simple algebra $A=\left( \frac{%
\begin{array}{c}
\\ 
\alpha ,\beta
\end{array}
}{%
\begin{array}{c}
K_{v},\varepsilon
\end{array}
}\right) $ is a division algebra if and only if $\beta \notin
N_{L_{w}/K_{v}}\left( L_{w}^{*}\right) .\Box \medskip $

From the above proposition and the Theorem 1.9, result that a $K_{v}-$cyclic
central simple algebra $A=\left( \frac{%
\begin{array}{c}
\\ 
\alpha ,\beta
\end{array}
}{%
\begin{array}{c}
K_{v},\varepsilon
\end{array}
}\right) $ is either  split or a division algebra.\medskip

In the following, we will study the central simple algebra $A=\left( \frac{%
\begin{array}{c}
\\ 
\alpha ,\,\,p^{3l}
\end{array}
}{%
\begin{array}{c}
K_{p},\varepsilon
\end{array}
}\right) ,$ where $p$ is a prime natural number, $p>3,l\in \mathbb{N}^{*},$ $\varepsilon $ is a primitive cubic root of unity, $K=%
\mathbb{Q}\left( \varepsilon \right) . $\newline
\medskip\newline
\textbf{Proposition 3.2.} \textit{Let }$p$\textit{\ be a prime natural
number,\thinspace \thinspace \thinspace }$p\equiv 2$ (\textit{mod} $3$)
\textit{\ and let be given the }$K_{p}-$\textit {\ algebra }
\textit{\ where }$l\in \mathbb{N}^{*},\alpha \in K,K=\mathbb{Q}\left(
\varepsilon \right) .$\textit{\ Let }$P$\textit{\ be a prime ideal of the
ring of integers of the field } $L=K\left(\sqrt[3]{\alpha}\right),$ \textit{\ lying above }$p.$\textit{\ Then }$p^{3l}$\textit{\ is a
norm from }$L_{P}^{*}$\textit{\ and the local Artin symbol }$\left( \frac{%
\begin{array}{c}
\\ 
L_{P}\,/\,K_{p}
\end{array}
}{%
\begin{array}{c}
(p^{3l)}
\end{array}
}\right) $\textit{\ is the identity.\medskip }

\textbf{Proof}. Since\thinspace $p\equiv 2$ (mod $3$), from Theorem
1.13., we obtain that $p$ is prime in the ring\ $\mathbb{Z}[\varepsilon ].$ It
results that cubic residual symbol $\left( \frac{%
\begin{array}{c}
\\ 
\alpha
\end{array}
}{%
\begin{array}{c}
p_{1}\mathbb{Z}[\varepsilon ]
\end{array}
}\right)_{3} =1,$ from\newline
\smallskip\newline
 Theorem 1.12, we have that $p$ is totally split in $\mathbb{A},$ where $%
\mathbb{A}$ is the ring of integers of the field $L=K\left(\sqrt[3]{\alpha}\right) :\,\,p\mathbb{A}=P_{1}P_{2}P_{3},\,\,P_{i}\in Spec\left( 
\mathbb{A}\right) ,i=\overline{1,3}.$\newline
We denote with $g$ the number of decomposition of the ideal $p\mathbb{A}$ in the extension $K\subset L.$
It results $g=3$ and knowing that $efg=[L:K]=3,$
then\thinspace \thinspace \thinspace $f=e=1.$ But $[L_{P}:K_{p}]=ef,$
therefore $L_{P}=K_{p},$ for each $P\in Spec\left( \mathbb{A}\right) ,P\mid p%
\mathbb{A}.$ In this case, we obtain that $p$ is the norm of \thinspace itself
\thinspace in \thinspace the \thinspace trivial extension of\thinspace $K_{p}
$ and the Artin symbol $\left( \frac{%
\begin{array}{c}
\\ 
L_{P}\,/\,K_{p}
\end{array}
}{%
\begin{array}{c}
(p^{3l})
\end{array}
}\right) $ is the identity.\medskip $\Box $

\textbf{Proposition 3.3.} \textit{Let }$p\,\,\,$\textit{be a prime natural
number,\thinspace \thinspace \thinspace }$p\equiv 1\,\,$(\textit{%
mod\thinspace \thinspace }$3$)\textit{\ and \thinspace let }$K_{p_{1}}-$%
\textit{\ algebra }$A=\left( \frac{%
\begin{array}{c}
\\ 
\alpha ,p^{3l}
\end{array}
}{%
\begin{array}{c}
K_{p_{1}},\varepsilon
\end{array}
}\right) ,$\textit{\ where }$l\in \mathbb{N}^{*},\alpha \in K,K=\mathbb{Q}\left(
\varepsilon \right) $\textit{\ and }$p_{1}$\textit{\ is a prime element in }$%
\mathbb{Z}[\varepsilon ],$\textit{\ }$p_{1}\mid p.$\textit{\ Let }$P$\textit{\
be a prime ideal in the ring of integers of the field } $L=K\left(\sqrt[3]{\alpha}\right),$ \textit{\ lying above }$p_{1}.$\textit{\ Then }$%
p^{3l}\in N_{L_{P}/K_{p_{1}}}\left( L_{P}^{*}\right) $\textit{\ and the
local Artin symbol \thinspace }$\left( \frac{%
\begin{array}{c}
\\ 
L_{P}\,/\,K_{p_{1}}
\end{array}
}{%
\begin{array}{c}
(p^{3l})
\end{array}
}\right) $\textit{\ is the identity in the Galois group }$Gal\left(
L_{P}/K_{p_{1}}\right) .\medskip $

\textbf{Proof.} From Theorem 1.13 and that $\mathbb{Z}[\varepsilon ]$ is a
principal ring, we have that the ideal $p\mathbb{Z}[\varepsilon ]=p_{1}\mathbb{Z}%
[\varepsilon ]\cdot p_{2}\mathbb{Z}[\varepsilon ],$ where $p_{1},p_{2}$ are
prime distinct elements in $\mathbb{Z}[\varepsilon ].$

We study the $K_{p_{1}}-$ algebra $A=\left( \frac{%
\begin{array}{c}
\\ 
\alpha ,p^{3l}
\end{array}
}{%
\begin{array}{c}
K_{p_{1}},\varepsilon
\end{array}
}\right) .$

\textbf{Case 1}. If the cubic residual symbol $\left( \frac{%
\begin{array}{c}
\\ 
\alpha
\end{array}
}{%
\begin{array}{c}
p_{1}\mathbb{Z}[\varepsilon ]
\end{array}
}\right)_{3} $ is a root of unity different\newline 
\smallskip\newline
from $1,$ from Theorem 1.12, we
obtain that the ideal $p_{1}\mathbb{A}\in Spec\left( \mathbb{A}\right) ,$where $%
\mathbb{A}$ is the ring of integers of the Kummer field $K\left(\sqrt[3]{\alpha}\right).$ So that $e=1,g=1$ and since $efg=[K\left(\sqrt[3]{\alpha}\right) :K]=3,$ it results that $f=3,$ who obviously divides $3l.$ From
Theorem 1.10, we obtain that $p^{3l}\in N_{L_{P}/K_{p_{1}}}\left(
L_{P}^{*}\right) .$ Using Theorem 1.11 and Proposition 3.1, we have
that the local Artin symbol \thinspace $\left( \frac{%
\begin{array}{c}
\\ 
L_{P}\,/\,K_{p_{1}}
\end{array}
}{%
\begin{array}{c}
(p^{3l})
\end{array}
}\right) $ is the identity in the Galois group $Gal\left(
L_{P}/K_{p_{1}}\right) $ and the algebra $A=\left( \frac{%
\begin{array}{c}
\\ 
\alpha ,p^{3l}
\end{array}
}{%
\begin{array}{c}
K_{p_{1}},\varepsilon
\end{array}
}\right) $ is not a division $K_{p_{1}}$ algebra.

\textbf{Case2.} If the cubic residual symbol $\left( \frac{%
\begin{array}{c}
\\ 
\alpha
\end{array}
}{%
\begin{array}{c}
p_{1}\mathbb{Z}[\varepsilon ]
\end{array}
}\right)_{3} =1,$ from Theorem 1.12,\newline 
\smallskip\newline
we obtain that $p_{1}\mathbb{A}%
=P_{1}P_{2}P_{3},\,\,P_{i}\in Spec\left( \mathbb{A}\right) ,i=\overline{1,3},$
therefore $g=3.$ But\thinspace $efg=[K\left(\sqrt[3]{\alpha}\right) :K]=3,$ therefore $e=f=1.$ Since $[L_{P}:K_{p_{1}}]=ef,$ we obtain
that $L_{P}=K_{p_{1}}$ for each $P\in Spec\left( \mathbb{A}\right) ,p\mid p_{1}%
\mathbb{A.}$ In this case, we have that $p_{1}$ is a norm of itself in the
trivial extension of $K_{p_{1}}$ and the local Artin symbol $\left( \frac{%
\begin{array}{c}
\\ 
L_{P}\,/\,K_{p_{1}}
\end{array}
}{%
\begin{array}{c}
(p^{3l})
\end{array}
}\right) $ is the identity.\medskip $\Box $%
\[
\]
\bigskip\\
\bigskip\\
\bigskip\\
  \textbf{Acknowledgements}\\
  The first author is indebted to Senior Research Fellow Tamas Szamuely, from Alfr$\acute{e}$d R$\acute{e}$nyi Institute of Mathematics from Budapest for the many helpful discussions.\\
  \bigskip\\
{\large References} 
\[
\]

[1] Acciaro, V., \textit{Solvability of Norm Equations over Cyclic Number
Fields of Prime Degree}, Mathematics of Computation, \textbf{65}(216)(1996),
1663-1674.\medskip

[2] Acciaro,V., Kluners, J., \textit{Computing Local Artin Maps, and
Solvability of Norm Equations}, Journal Symbolic Computation \textbf{11}%
(2000), 1-14.\medskip

[3] Elduque, E., \textit{Okubo algebras and twisted polynomials},
Contemporary Mathematics, \textbf{224}(1999), 101-109.\medskip

[4] Gille, P., Szamuely, T., \textit{Central Simple Algebras and Galois
Cohomology}, Cambridge University Press, 2006.\medskip

[5] Ireland,K., Rosen M. \textit{A Classical Introduction to Modern Number
Theory}, Springer Verlag, 1992.\medskip

[6] Janusz, G.J., \textit{Algebraic number fields}, Academic Press, London,
1973.\medskip

[7] Milne, J.S., \textit{Class Field Theory},
http://www.math.lsa.umich.edu/~jmilne.\medskip

[8] Pierce, R.S., \textit{Associative Algebras}, Springer Verlag,
1982.\medskip

[9] Schafer, R. D., \textit{An Introduction to Nonassociative Algebras},
Academic Press, New-York, 1966.
 
\end{document}